\documentclass{article}
\usepackage{amssymb,latexsym,amsfonts}

\title{Normalized Leonard pairs and Askey-Wilson relations}

\author{Raimundas Vid\=unas\thanks{
Supported by the 21 Century COE Programme "Development of Dynamic
Mathematics with High Functionality" of the Ministry of Education, Culture,
Sports, Science and Technology of Japan.}\\
\em Kyushu University}

\newtheorem{theorem}{Theorem}[section]
\newtheorem{lemma}[theorem]{Lemma}

\newtheorem{definition}[theorem]{Definition}

\newtheorem{remark}[theorem]{Remark}
\newtheorem{question}[theorem]{Question}
\newcommand{\proof}{{\bf Proof. }}
\newcommand{\qed}{\hfill $\Box$}
\newcommand{\equal}{&\!\!\!=\!\!\!&}

\newcommand{\ZZ}{{\Bbb Z}}

\newcommand{\fK}{{\Bbb K}}

\date{}

\begin{document}

\maketitle

\begin{abstract}
Let $V$ denote a vector space with finite positive dimension, and let
$(A,A^*)$ denote a Leonard pair on $V$. As is known, the linear
transformations $A$, $A^*$ satisfy the Askey-Wilson relations
\begin{eqnarray*}
A^2 A^*-\beta A A^*\!A+A^*\!A^2-\gamma\left( A A^*\!+\!A^*\!A
\right)-\varrho\,A^* \equal \gamma^*\!A^2+\omega A+\eta\,I,\\
A^*{}^2\!A-\beta A^*\!AA^*\!+AA^*{}^2\!-\gamma^*\!\left(A^*\!A\!+\!A
A^*\right)-\varrho^*\!A \equal \gamma A^*{}^2\!+\omega A^*\!+\eta^*I,
\end{eqnarray*}
for some scalars 
$\beta,\gamma,\gamma^*,\varrho,\varrho^*,\omega,\eta,\eta^*$. The scalar
sequence is unique if the dimension of $V$ is at least $4$.

If $c,c^*\!,t,t^*\!$ are scalars and $t,t^*\!$ are not zero, then
$(tA+c,t^*\!A^*\!+c^*)$ is a Leonard pair on $V$ as well. These affine
transformations can be used to bring the Leonard pair or its Askey-Wilson
relations into a convenient form. This paper presents convenient
normalizations of Leonard pairs by the affine transformations, and exhibits
explicit Askey-Wilson relations satisfied by them.\\ 

{\noindent\bf AMS 2000 MSC Classification}: 05E35, 33D45, 33C45.\\


{\noindent\bf E-mail address}: {\sf vidunas@math.kyushu-u.ac.jp}.\\
\end{abstract}

\newpage

\section{Introduction}

Throughout the paper, $\fK$ denotes an algebraically closed field. Apart
from one remark, we assume the characteristic of $\fK$ is not equal to 2.

Recall that a tridiagonal matrix is a square matrix which has nonzero
entries only on the main diagonal, on the superdiagonal and the subdiagonal.
A tridiagonal matrix is called irreducible whenever all entries on the
superdiagonal and subdiagonal are nonzero.

\begin{definition} \label{deflp} \rm
Let $V$ be a vector space over $\fK$ with finite positive dimension. By a
{\em Leonard pair} on $V$ we mean an ordered pair $(A,A^*)$, where $A:V\to
V$ and $A^*:V\to V$ are linear transformations which satisfy the following
two conditions:
\begin{enumerate}
\item[\it(i)] There exists a basis for $V$ with respect to which the matrix
representing $A^*$ is diagonal, and the matrix representing $A$ is
irreducible tridiagonal.%
\item[\it(ii)] There exists a basis for $V$ with respect to which the matrix
representing $A$ is diagonal, and the matrix representing $A^*$ is
irreducible tridiagonal.
\end{enumerate}
\end{definition}
\begin{remark} \rm
In this paper we do not use the conventional notation $A^*$ for the
conjugate-transpose of $A$. In a Leonard pair $(A,A^*)$, the linear
transformations $A$ and $A^*$ are arbitrary subject to the conditions {\it
(i)} and {\it (ii)} above.
\end{remark}

Leonard pairs occur in the theory of orthogonal polynomials, combinatorics,
the representation theory of the Lie algebra $sl_2$ or the quantum group
$U_q(sl_2)$. We refer to \cite{terwgen} 
as a survey on Leonard pairs, and as a source of further references.

We have the following result \cite[Theorem 1.5]{TerwVid}.
\begin{theorem} \label{lptheorem}
Let $V$ denote a vector space over $\fK$ with finite positive dimension. Let
$(A,A^*)$ be a Leonard pair on $V$. Then there exists a sequence of scalars
$\beta,\gamma,\gamma^*,\varrho,\varrho^*$, $\omega,\eta,\eta^*$ taken from
$\fK$ such that
\begin{eqnarray}  \label{askwil1}
A^2 A^*-\beta A A^*\!A+A^*\!A^2-\gamma\left( A A^*\!+\!A^*\!A
\right)-\varrho\,A^* \equal \gamma^*\!A^2+\omega A+\eta\,I,\\
\label{askwil2} A^*{}^2\!A-\beta A^*\!AA^*\!+AA^*{}^2-
\gamma^*\!\left(A^*\!A\!+\!A A^*\right)-\varrho^*\!A \equal \gamma
A^*{}^2+\omega A^*\!+\eta^*I.
\end{eqnarray}
The sequence is uniquely determined by the pair $(A,A^*)$ provided the
dimension of $V$ is at least $4$.
\end{theorem}
The equations (\ref{askwil1})--(\ref{askwil2}) are called the {\em
Askey-Wilson relations}. They first appeared in the work \cite{Zhidd} of
Zhedanov, where he showed that the Askey-Wilson polynomials give pairs of
infinite-dimensional matrices which satisfy the Askey-Wilson relations. We
denote this pair of equations by
$AW(\beta,\gamma,\gamma^*,\varrho,\varrho^*,\omega,\eta,\eta^*)$.

It is easy to notice that if $(A,A^*)$ is a Leonard pair, then
\begin{equation} \label{translation}
(tA+c,\,t^*\!A^*\!+c^*),\qquad \mbox{with}\ c,c^*,t,t^*\!\in\fK\mbox{ and }
t,t^*\neq 0,
\end{equation}
is a Leonard pair as well. We say that the two Leonard pairs 
are related by the {\em affine transformation} $(A,A^*)\mapsto
(tA+c,\,t^*\!A^*\!+c^*)$. 
Affine transformations act on Askey-Wilson relations as well, as explained
in Section \ref{normawrels} below. For example, if \mbox{$\beta\neq2$} then
the Askey-Wilson relations can be normalized so that $\gamma=0$ and
\mbox{$\gamma^*\!=0$}. Affine transformations can be used to normalize
Leonard pairs, parameter arrays representing them, or the Askey-Wilson
relations conveniently.

This paper presents convenient normalizations of Leonard pairs and their
Askey-Wilson relations. We generally assume that the dimension of the
underling vector space is at least 4, and use Terwilliger's classification
\cite{TerwLTparr} (or \cite[Section 35]{terwgen}) of parameter arrays
representing Leonard pairs. For parameter arrays of the $q$-type, we present
two normalizations: one that is close to Terwilliger's general expressions
in \cite{TerwLTparr}, and one where Askey-Wilson coefficients are normalized
most attractively. For other parameter arrays, we give one normalization.
This work is more of bookkeeping kind than of deep research. Examples of
Askey-Wilson relations for normalized Leonard pairs are given in
\cite{TerwVid}, \cite{hjalmarpaul}. Indirectly, Askey-Wilson relations for
Leonard pairs arising from certain distince regular graphs are computed in
\cite{Cur2hbipT}, \cite{Go}.

We note that Terwilliger's classification of parameter arrays by certain
families of orthogonal polynomials from the Askey-Wilson scheme can be
largely imitated to categorize Leonard pairs and Askey-Wilson relations; see
Sections \ref{lsparar} and \ref{classifics} below. We have the same types of
Leonard pairs and of Askey-Wilson relations, except that the quantum
$q$-Krawtchouk and affine $q$-Krawtchouk types are merged.

The paper is organized as follows. In the next section we discuss the
relation between Leonard pairs and parameter arrays. In Section
\ref{paawrel} we recall expressions of the Askey-Wilson coefficients in
(\ref{askwil1})--(\ref{askwil2}) in terms of parameter arrays. Section
\ref{normawrels} deals with possible normalizations of Askey-Wilson
relations. Sections \ref{qpararrays1} and \ref{qpararrays2} present two
normalizations of $q$-parameter arrays and Askey-Wilson relations for them.
Section \ref{otherpararrays} presents normalizations of other parameter
arrays and Askey-Wilson relations for them. In Section \ref{classifics} we
give a classification of Askey-Wilson relations consistent with the
classification of Leonard pairs. In Section \ref{conclusions} we discuss
uniqueness of normalizations of Leonard pairs and Askey-Wilson relations.

\section{Leonard pairs and parameter arrays}
\label{lsparar}

Leonard pairs are represented and classified by parameter arrays. More
precisely, parameter arrays are in one-to-one correspondence with {\em
Leonard systems} \cite[Definition 3.2]{terwgen}, and to each Leonard pair
one associates 4 Leonard systems or parameter arrays. 

From now on, let $d$ be a non-negative integer, and let $V$ be a vector
space with dimension $d+1$ over $\fK$.
\begin{definition} \rm \cite{hartwig}
Let $(A,A^*)$ denote a Leonard pair on $V$. Let $W$ denote a vector space
over $\fK$ with finite positive dimension, and let $(B,B^*)$ denote a
Leonard pair on $W$. By an {\em isomorphism of Leonard pairs} we mean an
isomorphism of vector spaces $\sigma:V\mapsto W$ such that $\sigma
A\sigma^{-1}=B$ and $\sigma A^*\sigma^{-1}=B^*$. We say that $(A,A^*)$ and
$(B,B^*)$ are {\em isomorphic} if there is an isomorphism of Leonard pairs
from $(A,A^*)$ to $(B,B^*)$.
\end{definition}

\begin{definition} \label{defpa} \rm
\cite{TerwLTparr} By a {\em parameter array} over $\fK$, of diameter $d$, we mean a sequence
\begin{equation} \label{paraar}
(\theta_0,\theta_1,\ldots,\theta_d;\;\theta_0^*,\theta_1^*,\ldots,\theta_d^*;\;
\varphi_1,\ldots,\varphi_d;\;\phi_1,\ldots,\phi_d)
\end{equation}
of scalars taken from $\fK$, that satisfy the following conditions:
\begin{enumerate}
\item[PA1.] $\theta_i\neq\theta_j$ and $\theta_i^*\neq\theta_j^*$ if $i\neq j$, for $0\le i,j\le d$.%
\item[PA2.] $\varphi_i\neq 0$ and $\phi_i\neq 0$, for $1\le i\le d$.%
\item[PA3.] $\displaystyle\varphi_i=\phi_1\sum_{j=0}^{i-1}
\frac{\theta_j-\theta_{d-j}}{\theta_0-\theta_d}+\left(\theta^*_i-\theta_0^*\right)
\left(\theta_{i-1}-\theta_d\right)$, for $1\le i\le d$.%
\item[PA4.] $\displaystyle\phi_i=\varphi_1\sum_{j=0}^{i-1}
\frac{\theta_j-\theta_{d-j}}{\theta_0-\theta_d}+\left(\theta^*_i-\theta_0^*\right)
\left(\theta_{d-i+1}\!-\theta_0\right)$, for $1\le i\le d$.%
\item[PA5.] The expressions
\[ \frac{\theta_{i-2}-\theta_{i+1}}{\theta_{i-1}-\theta_i},\qquad
\frac{\theta^*_{i-2}-\theta^*_{i+1}}{\theta^*_{i-1}-\theta^*_i}
\]
are equal and independent of $i$, for $2\le i\le d-1$.
\end{enumerate}
\end{definition}
\medskip

To get a Leonard pair from parameter array (\ref{paraar}), one
must choose a basis for $V$ and define the two linear transformations by the
following matrices (with respect to that basis):
\begin{equation} \label{split1}
\left(\begin{array}{ccccc} \theta_0 \\ 1 & \theta_1 \\ & 1 & \theta_2 \\
& & \ddots & \ddots \\ & & & 1 & \theta_d \end{array}\right),\qquad
\left(\begin{array}{ccccc} \theta^*_0 & \varphi_1 \\ & \theta^*_1 &
\varphi_2 \\ & & \theta^*_2 & \ddots \\
& & & \ddots  & \varphi_d \\ & & & & \theta^*_d \end{array}\right).
\end{equation}
Alternatively, the following two matrices define an isomorphic Leonard pair
on $V$:
\begin{equation} \label{split2}
\left(\begin{array}{ccccc} \theta_d \\ 1 & \theta_{d-1} \\ & 1 & \theta_{d-2} \\
& & \ddots & \ddots \\ & & & 1 & \theta_0 \end{array}\right),\qquad
\left(\begin{array}{ccccc} \theta^*_0 & \phi_1 \\ & \theta^*_1 &
\phi_2 \\ & & \theta^*_2 & \ddots \\
& & & \ddots  & \phi_d \\ & & & & \theta^*_d \end{array}\right).
\end{equation}

Conversely, if $(A,A^*)$ is a Leonard pair on $V$, there exists
\cite[Section 21]{terwgen} a basis for $V$ with respect to which the
matrices for $A$, $A^*$ have the bidiagonal forms in (\ref{split1}),
respectively. There exists another basis for $V$ with respect to which the
matrices for $A$, $A^*$ have the bidiagonal forms in (\ref{split2}),
respectively, with the same scalars $\theta_0,\theta_1,\ldots,\theta_d;
\theta^*_0,\theta^*_1,\ldots,\theta^*_d$. Then the following 4 sequences are
parameter arrays of diameter $d$:
\begin{eqnarray} \label{parray1}
(\theta_0,\theta_1,\ldots,\theta_d;\;\theta_0^*,\theta_1^*,\ldots,\theta_d^*;\;
\varphi_1,\ldots,\varphi_d;\;\phi_1,\ldots,\phi_d),\\
\label{parray2}
(\theta_0,\theta_1,\ldots,\theta_d;\;\theta_d^*,\ldots,\theta_{1}^*,\theta_0^*;\;
\phi_d,\ldots,\phi_1;\;\varphi_d,\ldots,\varphi_1),\\ \label{parray3}
(\theta_d,\ldots,\theta_{1},\theta_0;\;\theta_0^*,\theta_1^*,\ldots,\theta_d^*;\;
\phi_1,\ldots,\phi_d;\;\varphi_1,\ldots,\varphi_d),\\ \label{parray9}
(\theta_d,\ldots,\theta_{1},\theta_0;\;\theta_d^*,\ldots,\theta_{1}^*,\theta_0^*;\;
\varphi_d,\ldots,\varphi_1;\;\phi_d,\ldots,\phi_1).
\end{eqnarray}
If we apply to any of these parameter arrays the construction above, we get
back a Leonard pair isomorphic to $(A,A^*)$. These are all parameter arrays
which correspond to $(A,A^*)$ in this way.

The parameter arrays in (\ref{parray1})--(\ref{parray9}) are related by
permutations. The permutation group is isomorphic to $\ZZ_2\times\ZZ_2$. The
group action is without fixed points, since the eigenvalues $\theta_i$'s (or
$\theta_i^*$'s) are distinct. Let $\downarrow$ and $\Downarrow$ denote the
permutations which transform (\ref{parray1}) into, respectively,
(\ref{parray2}) and (\ref{parray3}). Observe that the composition
$\downarrow\Downarrow$ transforms (\ref{parray1}) into (\ref{parray9}). We
refer to the permutations $\downarrow$, $\Downarrow$ and
$\downarrow\Downarrow$ as {\em relation operators}, because in \cite[Section
4]{terwgen} the parameter arrays in (\ref{parray1})--(\ref{parray9})
corresponding to $(A,A^*)$ and the 4 similar parameter arrays corresponding
to the Leonard pair $(A^*,A)$ are called {\em relatives} of each other.

Parameter arrays are classified by Terwilliger in
\cite{TerwLTparr}; alternatively, see \cite[Section 35]{terwgen}.
For each parameter array, certain orthogonal polynomials
naturally occur in entries of the transformation matrix between two bases
characterized in Definition \ref{deflp} for the corresponding Leonard pair.
Terwilliger's classification largely mimics the terminating branch of
orthogonal polynomials in the Askey-Wilson scheme \cite{koekswart}.
Specifically, the classification comprises Racah, Hahn, Krawtchouk
polynomials and their $q$-versions, plus Bannai-Ito and orphan polynomials.
Classes of parameter arrays can be identified by the type of corresponding
orthogonal polynomials; we refer to them as {\em Askey-Wilson types}. The
type of a parameter array is unambiguously defined if $d\ge 3$. 
We recapitulate Terwilliger's classification in Sections \ref{qpararrays1}
through \ref{otherpararrays} by giving general normalized parameter arrays
of each type.

By inspecting Terwilliger's general parameter arrays \cite[Section
35]{terwgen}, one can observe that the relation operators $\downarrow$,
$\Downarrow$, $\downarrow\Downarrow$ do not change the Askey-Wilson type of
a parameter array (but only the free parameters such as $q,h,h^*,s$ there),
except that the $\Downarrow$ and $\downarrow\Downarrow$ relations mix up the
quantum $q$-Krawtchouk and affine $q$-Krawtchouk types. Consequently, given
a Leonard pair, all 4 associated parameter arrays have the same type, except
when parameter arrays of the quantum $q$-Krawtchouk or affine $q$-Krawtchouk
type occur. Therefore we can use the same classifying terminology for
Leonard pairs, except that we have to merge the quantum $q$-Krawtchouk and
affine $q$-Krawtchouk types. 

\section{Parameter arrays and AW relations}
\label{paawrel}

Let us consider a parameter array as in (\ref{parray1}). Suppose that the
corresponding Leonard pair satisfies Askey-Wilson relations
$AW(\beta,\gamma,\gamma^*,\varrho,\varrho^*,\omega,\eta,\eta^*)$. Note that
the Askey-Wilson relations are invariant under isomorphism of Leonard pairs.
Expressions for the 8 Askey-Wilson coefficients in terms of parameter arrays
are presented in \cite[Theorem 4.5 and Theorem 5.3]{TerwVid}. Here are the
formulas:
\begin{eqnarray} \label{betap1}
\beta+1\equal\frac{\theta_{i-2}-\theta_{i+1}}{\theta_{i-1}-\theta_i}=
\frac{\theta^*_{i-2}-\theta^*_{i+1}}{\theta^*_{i-1}-\theta^*_i},\\
\label{betap2}\gamma\equal\theta_{i-1}-\beta\theta_i+\theta_{i+1},\\
\label{betap3}\gamma^*\equal\theta^*_{i-1}-\beta\theta^*_i+\theta^*_{i+1},\\
\label{betap4}\varrho\equal\theta_i^2-\beta\,\theta_i\,\theta_{i-1}
+\theta_{i-1}^2-\gamma\,(\theta_i+\theta_{i-1}),\\
\label{betap5}\varrho^*\equal\theta^{*\,2}_i-\beta\theta^*_i\theta^*_{i-1}
+\theta^{*2}_{i-1}- \gamma^*(\theta^*_i+\theta^*_{i-1}),\\
\label{betap6}\omega\equal a_i\,(\theta^*_i-\theta^*_{i+1})+a_{i-1}\,
(\theta^*_{i-1}-\theta^*_{i-2})-\gamma\,(\theta^*_i+\theta^*_{i-1})\\
\label{betap7}\equal a^*_i\,(\theta_i-\theta_{i+1})+a^*_{i-1}\,
(\theta_{i-1}-\theta_{i-2})-\gamma^*\,(\theta_i+\theta_{i-1}),\\
\label{betap8}\eta\equal a^*_i\,(\theta_{i}\!-\!\theta_{i-1})\,
(\theta_{i}\!-\!\theta_{i+1})-\gamma^*\,\theta_i^2-\omega\,\theta_i,\\
\label{betap9}\eta^*\equal a_i\,(\theta^*_{i}\!-\!\theta^*_{i-1})\,
(\theta^*_{i}\!-\!\theta^*_{i+1})-\gamma\,\theta_i^*{}^2-\omega\,\theta^*_i.
\end{eqnarray}
The expressions for $\beta+1$ and $\omega$ are valid for $2\le i\le d-1$,
the expressions for $\varrho$, $\varrho^*$ are valid for $1\le i\le d$, and
the expressions for $\gamma,\gamma^*,\eta,\eta^*$ are valid for $1\le i\le
d-1$. The numbers $a_i,a_i^*$ 
are the diagonal entries in the tridiagonal forms of $A$, $A^*$ of Definition \ref{deflp};
see \cite[Section 7]{terwgen}. In terms of parameter arrays, we have
\cite[Section 10]{Terw24}:
\begin{eqnarray} \label{asipq0}
a_i\equal \theta_i+\frac{\varphi_i}{\theta^*_i-\theta^*_{i-1}}+
\frac{\varphi_{i+1}}{\theta^*_i-\theta^*_{i+1}}\\
\label{asipq1} \equal \theta_{d-i}+\frac{\phi_i}{\theta^*_i-\theta^*_{i-1}}+
\frac{\phi_{i+1}}{\theta^*_i-\theta^*_{i+1}},\\
\label{asipq2} a_i^*\equal
\theta_i^*+\frac{\varphi_i}{\theta_i-\theta_{i-1}}+
\frac{\varphi_{i+1}}{\theta_i-\theta_{i+1}}\\ \label{asipq}
\equal\theta^*_{d-i}+\frac{\phi_{d-i+1}}{\theta_i-\theta_{i-1}}+
\frac{\phi_{d-i}}{\theta_i-\theta_{i+1}}.
\end{eqnarray}
Here for $i\in\{0,d\}$ we should take
\begin{equation} \label{phiconv}
\varphi_0=0, \quad \varphi_{d+1}=0, \quad \phi_0=0, \quad\phi_{d+1}=0.
\end{equation}
The numbers $\theta_{-1},\theta_{d+1},\theta_{-1}^*,\theta_{d+1}$ can be
left undetermined. Surely, the Askey-Wilson coefficients are invariant under
the action of $\downarrow$, $\Downarrow$, $\downarrow\Downarrow$ on
parameter arrays.

As stated in Theorem \ref{lptheorem}, the coefficient sequence
$\beta,\gamma,\gamma^*,\varrho,\varrho^*,\omega,\eta,\eta^*$ is unique if
$d\ge 3$. If $d=2$, we can take $\beta$ freely and other coefficients get
determined uniquely. 
If $d=1$, we can take the 3 coefficients $\beta,\gamma,\gamma^*$ freely. If
$d=0$, we can take the 6 coefficients
$\beta,\gamma,\gamma^*,\varrho,\varrho^*,\omega$ freely.

\section{Normalized Askey-Wilson relations}
\label{normawrels}

Let $(A,A^*)$ denote a Leonard pair on $V$. Suppose that it satisfies the
Askey-Wilson relations
$AW(\beta,\gamma,\gamma^*,\varrho,\varrho^*,\omega,\eta,\eta^*)$. It can be
computed that Leonard pair (\ref{translation}) then satisfies
\begin{eqnarray} \label{newaw}
&&\hspace{-24pt}AW\left(
\beta,\;\gamma\,t+(2-\beta)\,c,\;\gamma^*t^*+(2-\beta)\,c^*,\;
\varrho\,t^2-2\gamma\,c\,t+(\beta-2)\,c^2,\right.\nonumber \\
&&\varrho^*t^*{}^2\!-2\gamma^*c^*t^*\!+(\beta\!-\!2)\,c^*{}^2,\,\omega\,t\,t^*\!
-2\gamma\,c^*t-2\gamma^*c\,t^*\!+2(\beta\!-\!2)\,c\,c^*,\nonumber\\
&&\eta\,t^2t^*-\varrho\,c^*t^2-\omega\,c\,t\,t^*
+\gamma^*c^2t^*+2\gamma^*c\,c^*t+(2-\beta)\,c^2c^*,\nonumber\\
&&\left.\eta^*t\,t^*{}^2\!-\varrho^*c\,t^*{}^2\!-\omega\,c^*t\,t^*
+\gamma\,c^*{}^2t+2\gamma\,c\,c^*t^*\!+(2-\beta)\,c\,c^*{}^2\right).
\end{eqnarray}
Note that $\beta$ stays invariant. The affine transformations
\begin{equation} \label{afftr}
(A,A^*)\mapsto (tA+c,\, t^*\!A^*\!+c^*),\qquad \mbox{with}\
c,c^*,t,t^*\!\in\fK,\ t,t^*\neq 0,
\end{equation}
can be used to normalize Leonard pairs so that their Askey-Wilson relations
would have a simple form. We refer to a transformations of the form
$(A,A^*)\mapsto (A+c,A^*+c^*)$ as an {\em affine translation}, and to a
transformation of the form $(A,A^*)\mapsto (tA,\,t^*A^*)$ as an {\em affine
scaling}. Generally, we can use an affine translation to set some two
Askey-Wilson coefficients to zero, and then use an affine scaling to
normalize some two nonzero coefficients. Specifically, by affine
translations we can achieve the following.
\begin{lemma} \label{normrules}
The Askey-Wilson relations
$AW(\beta,\gamma,\gamma^*,\varrho,\varrho^*,\omega,\eta,\eta^*)$ can be
normalized as follows:
\begin{enumerate}
\item If $\beta\neq 2$, we can set $\gamma=0$, $\gamma^*=0$.%
\item If $\beta=2$, $\gamma\neq 0$, $\gamma^*\!\neq 0$, we can set $\varrho=0$, $\varrho^*=0$.%
\item If $\beta=2$, $\gamma=0$, $\gamma^*\!\neq 0$, we can set $\varrho^*=0$, $\omega=0$.%
\item If $\beta=2$, $\gamma^*\!=0$, $\gamma\neq 0$, we can set $\varrho=0$, $\omega=0$.%
\item If $\beta=2$, $\gamma=0$, $\gamma^*=0$, $\omega^2\neq\varrho\varrho^*$, we can set $\eta=0$, $\eta^*=0$.%
\item If $\beta=2$, $\gamma=0$, $\gamma^*=0$, $\mbox{\rm
rk}{\omega\;\; \varrho\;\;\eta\choose\,\varrho^*\,\;\omega\;\eta^*}\le 1$, we can set $\eta=0$, $\eta^*=0$.%
\item Otherwise, we have
\[
\beta=2,\quad \gamma=0,\quad \gamma^*=0,\quad
\omega^2=\varrho\varrho^*,\quad \mbox{\rm rk}{\omega\quad
\varrho\quad\eta\,\choose\varrho^*\quad\!\!\omega\quad\eta^*}=2.
\]
Then we can set either $\eta=0$ or $\eta^*=0$, but not both.
\end{enumerate}
In the first $5$ cases, there is a unique affine translation to make the
normalization. In the last $2$ cases, there are infinitely many
normalizations by affine translations.
\end{lemma}
\proof The first 4 cases are straightforward, including the uniqueness
statement. If $\beta=2$, $\gamma=0$, $\gamma^*=0$, the new Askey-Wilson
relations (\ref{newaw}) are
\begin{eqnarray*} \textstyle
AW\left( 2,\,0,\,0,\,\varrho\,t^2,\,\varrho^*t^*{}^2,\,
\omega\,tt^*\!,\,\big(\eta-\omega\,a-\varrho\,a^*\big)\,t^2t^*,\,
\big(\eta^*-\varrho^*a-\omega\,a^*\big)\,t\,t^*{}^2\right),
\end{eqnarray*}
where $a=c/t$ and $a^*=c^*/t^*$. To set the last two parameters to zero, we
have to solve two linear equations in $a,a^*$. If we have $\det{\omega\;\;
\varrho\choose\,\varrho^*\; \omega}\neq 0$, the solution is unique.
Otherwise we have either infinitely many or none solutions, which leads us
to the last two cases.
\qed\\

As it turns out, cases 6 and 7 of Lemma \ref{normrules} do not occur for
Askey-Wilson relations satisfied by Leonard pairs if $d\ge 3$. See part 3 of
Theorem \ref{classth} below.

In Section \ref{qpararrays1}, we normalize the general $q$-parameter arrays
in Terwilliger's classification \cite[Section 35]{terwgen} with most handy
changes in the explicit expressions. We use the following simplest action of
(\ref{afftr}) on parameter arrays, consistent with the transformation of
Leonard pairs:
\begin{equation} \label{scalingtr}
\theta_i\mapsto t\,\theta_i+c,\qquad \theta^*_i\mapsto t^*\theta^*_i+c^*,
\qquad \varphi_i\mapsto t\,t^*\varphi_i, \qquad \phi_i\mapsto t\,t^*\phi_i.
\end{equation}
It turns out that the corresponding Askey-Wilson relations follow the
specification of part 1 of Lemma \ref{normrules} immediately.

Suppose that we normalized a pair of Askey-Wilson relations to satisfy
implications of Lema \ref{normrules}, and suppose that cases 6 and 7 do not
apply. Then the only affine transformations which preserve two specified
zero coefficients are affine scalings. One can use affine scalings to
normalize some two nonzero coefficients to convenient values. Sections
\ref{qpararrays2} and \ref{otherpararrays} present such normalized parameter
arrays that in their Askey-Wilson relations two nonzero coefficients are
basically constants. (More precisely, in the $q$-cases they depend on $q$,
or equivalently, on $\beta$.) The scaling normalization is explained more
thoroughly in Section \ref{classifics}.

\section{Normalized $q$-parameter arrays}
\label{qpararrays1}

Here we present the most straightforward normalizations of the general
parameter arrays in \cite[Section 35]{terwgen} with the $q$-parameter. Lemma
\ref{awnormals} below gives the Askey-Wilson relations for the corresponding
Leonard pairs. The Askey-Wilson relations turn out to be normalized
according to part 1 of Lemma \ref{normrules}.
\begin{lemma} \label{normlps}
The parameter arrays in {\rm\cite[Examples 35.2--35.8]{terwgen}} can be
normalized by affine transformations $(\ref{scalingtr})$ to the following
forms:
\begin{itemize}
\item The $q$-Racah case: $\displaystyle\theta_i=q^{-i}+s\,q^{i+1},\quad
\theta^*_i=q^{-i}+s^*q^{i+1}$.
\begin{eqnarray*}
\varphi_i\equal q^{1-2i}\left(1-q^i\right)\left(1-q^{i-d-1}\right)
\left(1-r\,q^i\right)\left(r-s\,s^*q^{d+1+i}\right)\big/r,\\
\phi_i\equal q^{1-2i}\left(1-q^i\right)\nonumber
\left(1-q^{i-d-1}\right)\left(r-s^*q^i\right)\left(s\,q^{d+1}-r\,q^i\right)\big/r.
\end{eqnarray*}
\item The $q$-Hahn case: $\displaystyle\theta_i=q^{-i},\quad
\theta^*_i=q^{-i}+s^*q^{i+1}$,
\begin{eqnarray*}
\varphi_i\equal q^{1-2i}\,\left(1-q^i\right)\left(1-q^{i-d-1}\right)\left(1-r\,q^i\right),\\
\phi_i\equal-q^{1-i}\,\left(1-q^i\right)\left(1-q^{i-d-1}\right)\left(r-s^*q^i\right).
\end{eqnarray*}
\item The dual $q$-Hahn case: $\displaystyle
\theta_i=q^{-i}+s\,q^{i+1},\quad\theta^*_i=q^{-i}$,
\begin{eqnarray*}
\varphi_i\equal q^{1-2i}\,\left(1-q^i\right)\left(1-q^{i-d-1}\right)\left(1-r\,q^i\right),\\
\phi_i\equal
q^{d+2-2i}\,\left(1-q^i\right)\left(1-q^{i-d-1}\right)\left(s-r\,q^{i-d-1}\right).
\end{eqnarray*}
\item The $q$-Krawtchouk case: $\displaystyle\theta_i=q^{-i},\quad
\theta^*_i=q^{-i}+s^*q^{i+1}$,
\begin{eqnarray*}
\varphi_i\equal q^{1-2i}\,\left(1-q^i\right)\left(1-q^{i-d-1}\right),\\
\phi_i\equal s^*q\left(1-q^i\right)\left(1-q^{i-d-1}\right).
\end{eqnarray*}
\item The dual $q$-Krawtchouk case: $\displaystyle
\theta_i=q^{-i}+s\,q^{i+1},\quad\theta^*_i=q^{-i}$,
\begin{eqnarray*}
\varphi_i\equal q^{1-2i}\,\left(1-q^i\right)\left(1-q^{i-d-1}\right),\\
\phi_i\equal s\,q^{d+2-2i}\,\left(1-q^i\right)\left(1-q^{i-d-1}\right).
\end{eqnarray*}
\item The quantum $q$-Krawtchouk case: $\theta_i=q^{i+1}$,\,
$\theta^*_i=q^{-i}$,
\begin{eqnarray*}
\varphi_i\equal -r\,q^{1-i}\,\left(1-q^i\right)\left(1-q^{i-d-1}\right),\\
\phi_i\equal q^{d+2-2i}\,
\left(1-q^i\right)\left(1-q^{i-d-1}\right)\left(1-r\,q^{i-d-1}\right).
\end{eqnarray*}
\item The affine $q$-Krawtchouk case: $\theta_i=q^{-i}$,\,
$\theta^*_i=q^{-i}$,
\begin{eqnarray*}
\varphi_i\equal q^{1-2i}\,\left(1-q^i\right)\left(1-q^{i-d-1}\right)\left(1-r\,q^i\right),\\
\phi_i\equal-r\,q^{1-i}\,\left(1-q^i\right)\left(1-q^{i-d-1}\right).
\end{eqnarray*}
\end{itemize}
In each case, $q,s,s^*,r$ are nonzero scalar parameters such that
$\theta_i\neq \theta_j$, $\theta^*_i\neq\theta^*_j$ for $0\le i<j\le d$, and
$\varphi_i\neq 0$, $\phi_i\neq 0$ for $1\le i \le d$.
\end{lemma}
\proof By affine translations, we adjust Terwilliger's parameters
$\theta_0,\theta^*_0$ so that we have only summands depending on $i$ in the
expanded expressions for $\theta_i$, $\theta^*_i$ in \cite[Examples
35.2--35.8]{terwgen}. By affine scalings, we set Terwilliger's parameters
$h,h^*$ to the value 1. In the quantum $q$-Krawtchouk case \cite[Example
35.5]{terwgen} there is no parameter $h$, so we set $s=1$. Other parameters
remain unchanged, except that in the $q$-Racah case we rename $r_1$ to $r$
and set $r_2=s\,s^*q^{d+1}/r$.\qed

\begin{lemma} \label{awnormals}
Let $q,s,s^*,r$ denote the same scalar parameters as in the previous lemma.
We use the following notations:
\begin{eqnarray}\label{normconsts}
S=s\,q^{d+1}+1, \qquad S^*=s^*q^{d+1}+1, \qquad
R=r+\frac{s\,s^*\,q^{d+1}}r,\\ 
Q=q^{d+1}+1,\qquad K=-\frac{(q^2\!-1)^2}q, \quad
K^*=\frac{(q-1)^2}{q^{d+1}}.
\end{eqnarray}
The Askey-Wilson relations for the parameter arrays of Lemma $\ref{normlps}$
are: 
\begin{itemize}
\item For the $q$-Racah case:
\begin{eqnarray}
&&\hspace{-24pt}AW\big(q+q^{-1},\,0,\,0,\,s\,K,\,s^*K,\,-K^*\!\left(S\,S^*\!+R\,Q\right),\nonumber\\
&&(q+1)K^*\!\left( S\,R+s\,S^*Q\right),\,
(q+1)K^*\!\left(S^*R+s^*S\,Q\right)\big).
\end{eqnarray}
\item For the $q$-Hahn case:
\begin{eqnarray}
&&\hspace{-24pt}AW\big(q+q^{-1},\,0,\,0,\,0,\,s^*K,\,-K^*\!\left(S^*+r\,Q\right),\nonumber\\
&&(q+1)K^*r,\, (q+1)K^*\!\left(S^*r+s^*Q\right)\big).
\end{eqnarray}
\item For the dual $q$-Hahn case:
\begin{eqnarray}
&&\hspace{-24pt}AW\big(q+q^{-1},\,0,\,0,\,s\,K,\,0,\,-K^*\!\left(S+r\,Q\right),\nonumber\\
&&\qquad (q\!+\!1)K^*\!\left(S\,r+s\,Q\right),\,(q\!+\!1)K^*r\big).
\end{eqnarray}
\item For the $q$-Krawtchouk case:
\begin{eqnarray}
AW\big(q+q^{-1},\,0,\,0,\,0,\,s^*K,\,-K^*\!S^*,\,0,\,(q\!+\!1)K^*s^*Q\big).
\end{eqnarray}
\item For the dual $q$-Krawtchouk case:
\begin{eqnarray}
AW\big(q+q^{-1},\,0,\,0,\,s\,K,\,0,\,-K^*\!S,\,(q\!+\!1)K^*s\,Q,\,0\big).
\end{eqnarray}
\item For the quantum $q$-Krawtchouk case:
\begin{eqnarray}
AW\big(q+q^{-1},\,0,\,0,\,0,\,0,\,-K^*\!\left(q^{d+1}\!+r\,Q\right),\,
(q\!+\!1)(q\!-\!1)^2r,\,(q\!+\!1)K^*r\big).
\end{eqnarray}
\item For the affine $q$-Krawtchouk case:
\begin{eqnarray}
AW\big(q+q^{-1},\,0,\,0,\,0,\,0,\,-K^*\!\left(1+r\,Q\right),\,
(q\!+\!1)K^*r,\,(q\!+\!1)K^*r\big).\quad
\end{eqnarray}
\end{itemize}
\end{lemma}
\proof Direct computations with formulas (\ref{betap1})--(\ref{asipq2}).\qed

\section{Alternative normalized $q$-arrays}
\label{qpararrays2}

Here we present alternative normalizations of the general parameter arrays
in \cite[Section 35]{terwgen} with the general $q$-parameter. The parameters
are rescaled, and the free parameters $q,s,s^*,r$ are different. In
particular, the $q$ of the previous section is replaced by $q^2$. The
normalization for the $q$-Racah case is proposed in \cite{hjalmarpaul}.

The corresponding Askey-Wilson relations are normalized according to part 1
of Lemma \ref{normrules}. Advantages of this normalization are: the two
nonzero values normalized by affine scaling are $q$-constants;  
expressions for normalized parameter arrays are more symmetric;
the set of normalized parameter arrays is preserved by the $\downarrow$,
$\Downarrow$, $\downarrow\Downarrow$ operations (see Section
\ref{conclusions}).
\begin{lemma} \label{normlps1}
The parameter arrays in {\rm\cite[Examples 35.2--35.8]{terwgen}} can be
normalized by affine transformations $(\ref{scalingtr})$ to the following
forms:
\begin{itemize}
\item The $q$-Racah case:
$\displaystyle\theta_i=s\,q^{d-2i}+\frac{q^{2i-d}}{s},\quad
\theta^*_i=s^*q^{d-2i}+\frac{q^{2i-d}}{s^*}$.
\begin{eqnarray*}
\varphi_i\equal\frac{q^{2d+2-4i}}{s\,s^*r}\left(1-q^{2i}\right)\left(1-q^{2i-2d-2}\right)
\left(s\,s^*-r\,q^{2i-d-1}\right)\left(s\,s^*r-q^{2i-d-1}\right),\\
\phi_i\equal\frac{q^{2d+2-4i}}{s\,s^*r}\left(1-q^{2i}\right)\nonumber
\left(1-q^{2i-2d-2}\right)\left(s^*r-s\,q^{2i-d-1}\right)\left(s^*-s\,r\,q^{2i-d-1}\right).
\end{eqnarray*}
\item The $q$-Hahn case: $\displaystyle\theta_i=r\,q^{d-2i},\quad
\theta^*_i=s^*q^{d-2i}+\frac{q^{2i-d}}{s^*}$,
\begin{eqnarray*}
\varphi_i\equal \frac{q^{2d+2-4i}}{r}\,
\left(1-q^{2i}\right)\left(1-q^{2i-2d-2}\right)\left(s^*r^2-q^{2i-d-1}\right),\\
\phi_i\equal-\frac{q^{d+1-2i}}{r\,s^*}\,
\left(1-q^{2i}\right)\left(1-q^{2i-2d-2}\right)\left(s^*-r^2q^{2i-d-1}\right).
\end{eqnarray*}
\item The dual $q$-Hahn case: $\displaystyle
\theta_i=s\,q^{d-2i}+\frac{q^{2i-d}}{s},\quad\theta^*_i=r\,q^{d-2i}$,
\begin{eqnarray*}
\varphi_i\equal \frac{q^{2d+2-4i}}{r}\,
\left(1-q^{2i}\right)\left(1-q^{2i-2d-2}\right)\left(s\,r^2-q^{2i-d-1}\right),\\
\phi_i\equal
\frac{q^{2d+2-4i}}{r\,s}\,\left(1-q^{2i}\right)\left(1-q^{2i-2d-2}\right)\left(r^2-s\,q^{2i-d-1}\right).
\end{eqnarray*}
\item The $q$-Krawtchouk: $\displaystyle\theta_i=q^{d-2i},\quad
\theta^*_i=s^*q^{d-2i}+\frac{q^{2i-d}}{s^*}$,
\begin{eqnarray*}
\varphi_i\equal s^*\,q^{2d+2-4i}\,\left(1-q^{2i}\right)\left(1-q^{2i-2d-2}\right),\\
\phi_i\equal\frac{1}{s^*}\,\left(1-q^{2i}\right)\left(1-q^{2i-2d-2}\right).
\end{eqnarray*}
\item The dual $q$-Krawtchouk: $\displaystyle
\theta_i=s\,q^{d-2i}+\frac{q^{2i-d}}{s},\quad\theta^*_i=q^{d-2i}$,
\begin{eqnarray*}
\varphi_i\equal s\;q^{2d+2-4i}\,\left(1-q^{2i}\right)\left(1-q^{2i-2d-2}\right),\\
\phi_i\equal\frac{q^{2d+2-4i}}{s}\,\left(1-q^{2i}\right)\left(1-q^{2i-2d-2}\right).
\end{eqnarray*}
\item The quantum $q$-Krawtchouk: $\theta_i=r\,q^{2i-d}$,\,
$\theta^*_i=r\,q^{d-2i}$,
\begin{eqnarray*}
\varphi_i\equal -\frac{q^{d+1-2i}}{r}\,\left(1-q^{2i}\right)\left(1-q^{2i-2d-2}\right),\\
\phi_i\equal \frac{q^{2d+2-4i}}{r}\,
\left(1-q^{2i}\right)\left(1-q^{2i-2d-2}\right)\left(r^3-q^{2i-d-1}\right).
\end{eqnarray*}
\item The affine $q$-Krawtchouk: $\theta_i=r\,q^{d-2i}$,\,
$\theta^*_i=r\,q^{d-2i}$,
\begin{eqnarray*}
\varphi_i\equal \frac{q^{2d+2-4i}}r\,\left(1-q^{2i}\right)\left(1-q^{2i-2d-2}\right)\left(r^3-q^{2i-d-1}\right),\\
\phi_i\equal-\frac{q^{d+1-2i}}r\,\left(1-q^{2i}\right)\left(1-q^{2i-2d-2}\right).
\end{eqnarray*}
\end{itemize}
In each case, $q,s,s^*,r$ are nonzero scalar parameters such that
$\theta_i\neq \theta_j$, $\theta^*_i\neq\theta^*_j$ for $0\le i<j\le d$, and
$\varphi_i\neq 0$, $\phi_i\neq 0$ for $1\le i \le d$.
\end{lemma}
\proof In every case of Lemma \ref{normlps}, we substitute
\[
q\mapsto q^2,\qquad s\mapsto\frac1{s^2\,q^{2d+2}},\qquad
s^*\mapsto\frac1{s^*{}^2q^{2d+2}}.
\]
Besides, in the $q$-Racah, $q$-Hahn, dual $q$-Hahn, quantum $q$-Krawtchouk
and affine $q$-Krawtchouk cases we substitute $r$ by, respectively,
\[
\frac{r}{s\,s^*q^{d+1}},\quad \frac{1}{s^*r^2q^{d+1}},\quad
\frac{1}{s\,r^2q^{d+1}},\quad \frac{q^{d+1}}{r^3},\quad
\frac{1}{r^3q^{d+1}}.
\]
After that, we apply affine scaling. We use formula (\ref{scalingtr}) with
$c=0,c^*=0$ and $(t,t^*)$ equal to, respectively in the listed order,
\begin{eqnarray*}
(s\,q^d,s^*q^d), \quad (r\,q^d,s^*q^d), \quad (s\,q^d,r\,q^d), \quad
(q^d,s^*q^d), \quad (s\,q^d,q^d),\\
(r\,q^{-d-2},r\,q^d), \quad (r\,q^d,r\,q^d).
\end{eqnarray*}
\qed

\begin{lemma} \label{awnormals1}
As in the previous lemma, let $q,s,s^*,r$ denote nonzero scalar parameters.
We use the following notations:
\begin{eqnarray}\label{normconsts1}
Q_j=q^j+q^{-j},\qquad Q_j^*=q^j-q^{-j}, \quad \mbox{for $j=1,2,\ldots$},\\
S=s+\frac{1}s, \qquad S^*=s^*+\frac{1}{s^*}, \qquad R=r+\frac{1}{r}.
\end{eqnarray}
The Askey-Wilson relations for the parameter arrays of Lemma
$\ref{normlps1}$ are: 
\begin{itemize}
\item For the $q$-Racah case:
\begin{eqnarray}
&&\hspace{-24pt}AW\big(Q_2,\,
0,\,0,\,-Q_2^*{}^2,\,-Q_2^*{}^2,
\,-Q_1^*{}^2\!\left(S\,S^*\!+Q_{d+1}R\right),\nonumber\\
&&Q_1Q_1^*{}^2\!\left( S\,R+Q_{d+1}S^*\right),\,
Q_1Q_1^*{}^2\!\left(S^*R+Q_{d+1}S\right)\big).
\end{eqnarray}
\item For the $q$-Hahn case:
\begin{eqnarray}
&&\hspace{-24pt}AW\big(Q_2,\,0,\,0,\,0,\,-Q^*_2{}^2,\,-Q_1^*{}^2\!\left(S^*r+Q_{d+1}r^{-1}\right),\nonumber\\
&&Q_1Q_1^*{}^2,\, Q_1Q_1^*{}^2\!\left(S^*r^{-1}+Q_{d+1}r\right)\big).
\end{eqnarray}
\item For the dual $q$-Hahn case:
\begin{eqnarray}
&&\hspace{-24pt}AW\big(Q_2,\,0,\,0,\,-Q^*_2{}^2,\,0,\,-Q_1^*{}^2\!\left(S\,r+Q_{d+1}r^{-1}\right),
\nonumber\\ &&\qquad
Q_1Q_1^*{}^2\!\left(S\,r^{-1}\!+Q_{d+1}r\right),\,Q_1Q_1^*{}^2\big).
\end{eqnarray}
\item For the $q$-Krawtchouk case:
\begin{eqnarray}
AW\big(Q_2,\,0,\,0,\,0,\,-Q_2^*{}^2,\,-Q_1^*{}^2S^*,\,0,\,Q_1Q_1^*{}^2Q_{d+1}\big).
\end{eqnarray}
\item For the dual $q$-Krawtchouk case:
\begin{eqnarray}
AW\big(Q_2,\,0,\,0,\,-Q_2^*{}^2,\,0,\,-Q_1^*{}^2S,\,Q_1Q_1^*{}^2Q_{d+1},\,0\big).
\end{eqnarray}
\item For the quantum $q$-Krawtchouk and affine $q$-Krawtchouk cases:
\begin{eqnarray}
AW\big(Q_2,\,0,\,0,\,0,\,0,\,-Q_1^*{}^2\!\left(r^2+Q_{d+1}r^{-1}\right),\,
Q_1Q_1^*{}^2,\,Q_1Q_1^*{}^2\big).\quad
\end{eqnarray}
\end{itemize}
\end{lemma}
\proof Transform the Askey-Wilson relations in Lemma \ref{awnormals} with
the same substitutions and affine scalings as in the previous proof. In the notation
of this lemma, the expressions $S,S^*,R,Q,K,K^*$ of
Lemma \ref{awnormals} get replaced by, respectively, $S/s$, $S^*/s^*$,
$R/q^{d+1}s\,s^*$, $q^{d+1}Q_{d+1}$, $-q^2Q_2^*{}^2$ and $q^{-2d}Q_1^*{}^2$.
\qed

\section{Other parameter arrays}
\label{otherpararrays}

Here we present normalizations of the remaining general parameter arrays in
\cite[Section 35]{terwgen}. The corresponding Askey-Wilson relations are
normalized according to Lemma \ref{normrules}, and two nonzero values are
constants. Since we generally assume that char $\fK\neq 2$, the orphan case
is missing in the lemmas below. It is briefly discussed in Remark \ref{orphan}.
\begin{lemma} \label{normlps2}
The parameter arrays in {\rm\cite[Examples 35.9--35.13]{terwgen}} can be
normalized by affine transformations $(\ref{scalingtr})$ to the following
forms:
\begin{itemize}
\item The Racah case: $\theta_i=(i+u)(i+u+1)$,
$\theta^*_i=(i+u^*)(i+u^*+1)$,
\begin{eqnarray*}
\varphi_i\equal i\,(i-d-1)\,(i+u+u^*-v)\,(i+u+u^*+d+1+v),\\
\phi_i\equal i\,(i-d-1)\,(i-u+u^*+v)\,(i-u+u^*-d-1-v).
\end{eqnarray*}
\item The Hahn case: $\theta_i=i+v-\frac{d}2$,
$\theta^*_i=(i+u^*)(i+u^*+1)$,
\begin{eqnarray*}
\varphi_i\equal i\,(i-d-1)\,(i+u^*\!+2v),\\
\phi_i\equal -i\,(i-d-1)\,(i+u^*\!-2v).
\end{eqnarray*}
\item The dual Hahn case: $\theta_i=(i+u)(i+u+1)$,
$\theta^*_i=i+v-\frac{d}2$,
\begin{eqnarray*}
\varphi_i\equal i\,(i-d-1)\,(i+u+2v),\\
\phi_i\equal i\,(i-d-1)\,(i-u+2v-d-1).
\end{eqnarray*}
\item The Krawtchouk case: $\theta_i=i-\frac{d}2$, $\theta^*_i=i-\frac{d}2$,
\begin{eqnarray*}
\varphi_i\equal v\,i\,(i-d-1),\\
\phi_i\equal (v-1)\,i\,(i-d-1).
\end{eqnarray*}
\item The Bannai-Ito case: $\theta_i=(-1)^i\left(i+u-\frac{d}2\right)$,
$\theta^*_i=(-1)^i\left(i+u^*\!-\frac{d}2\right)$,
\begin{eqnarray*}
\varphi_i=\left\{ \begin{array}{cl} -i\left(i+u+u^*\!+v-\frac{d+1}2\right),&
\mbox{for $i$ even, $d$ even}.\\
-(i-d-1)\left(i+u+u^*\!-v-\frac{d+1}2\right), &
\mbox{for $i$ odd, $d$ even}.\\
-i\,(i-d-1), & \mbox{for $i$ even, $d$ odd}.\\
v^2-\left(i+u+u^*\!-\frac{d+1}2\right)^2, &
 \mbox{for $i$ odd, $d$ odd}.\end{array} \right.
\end{eqnarray*}
\begin{eqnarray*}
\phi_i=\left\{ \begin{array}{cl} i\left(i-u+u^*\!-v-\frac{d+1}2\right),& \mbox{for $i$ even, $d$ even}.\\
(i-d-1)\left(i-u+u^*\!+v-\frac{d+1}2\right), & \mbox{for $i$ odd, $d$ even}.\\
-i\,(i-d-1), & \mbox{for $i$ even, $d$ odd}.\\
v^2-\left(i-u+u^*\!-\frac{d+1}2\right)^2, &
 \mbox{for $i$ odd, $d$ odd}.\end{array} \right.
\end{eqnarray*}
\end{itemize}
In each case, $u,u^*,v$ are scalar parameters such that $\theta_i\neq
\theta_j$, $\theta^*_i\neq\theta^*_j$ for $0\le i<j\le d$, and $\varphi_i\neq
0$, $\phi_i\neq 0$ for $1\le i \le d$.
\end{lemma}
\proof Like in the proof of Lemma \ref{normlps}, we adjust Terwilliger's
parameters $\theta_0,\theta^*_0$ by an affine translation, and then adjust
other two parameters by an affine scaling. We also make linear substitutions for
the remaining parameters. In the Racah case, we substitute
\[
s\mapsto 2u,\quad s^*\!\mapsto 2u^*,\quad r_1\mapsto u+u^*-v, \quad \mbox{so
that} \quad r_2=u+u^*+d+1+v.
\]
Then we adjust $\theta_0=u^2+u$, $\theta^*_0=u^*{}^2\!+u^*\!$, $h=1$,
$h^*\!=1$. In the Hahn case, we substitute $s^*\!\mapsto 2u^*$, $r\mapsto
u^*\!+2v$ and adjust $\theta_0=v-\frac{d}2$, $\theta^*_0=u^*{}^2\!+u^*\!$,
$h^*\!=1$, $s=1$. In the dual Hahn case, we substitute $s\mapsto 2u$,
$r\mapsto u+2v$ and adjust $\theta_0=u^2+u$, $\theta^*_0=v-\frac{d}2$,
$h=1$, $s=1$. In the Krawtchouk case, we substitute $r\mapsto v$ and adjust
$\theta_0=-\frac{d}2$, $\theta^*_0=-\frac{d}2$, $s=1$, $s^*\!=1$. In the
Bannai-Ito case, we substitute
\[ \textstyle
s\mapsto d+1-2u,\quad s^*\!\mapsto d+1-2u^*,\quad r_1\mapsto
u+u^*+v-\frac{d+1}2, 
\]
so that $r_2\mapsto u+u^*-v-\frac{d+1}2$, and adjust $\theta_0=u-\frac{d}2$,
$\theta^*_0=u^*\!-\frac{d}2$, $h=\frac12$, $h^*\!=\frac12$. \qed

\begin{lemma} \label{awnormals2}
Let $u,u^*,v$ denote the same scalar parameters as in the previous lemma.
The Askey-Wilson relations for the parameter arrays of Lemma
$\ref{normlps2}$
are: 
\begin{itemize}
\item For the Racah case:
\begin{eqnarray}
&&\hspace{-24pt}
AW\big(2,\,2,\,2,\,0,\,0,-2u^2-2u^*{}^2\!-2v^2-2(d\!+\!1)(u+u^*+v)-2d^2-4d,\nonumber\\
&&2\,u\left(u+d+1\right)\left(v-u^*\right)\left(v+u^*+d+1\right), \\
&&2\,u^*\!\left(u^*\!+d+1\right)\left(v-u\right)\left(v+u+d+1\right)\big).\nonumber
\end{eqnarray}
\item For the Hahn case:
\begin{eqnarray}
\textstyle
AW\big(2,\,0,\,2,\,1,\,0,\,0,\,-(u^*\!+1)(u^*\!+d)-2v^2-\frac{d^2}2,\,-4u^*(u^*\!+d+1)v\big).
\end{eqnarray}
\item For the dual Hahn case:
\begin{eqnarray}
\textstyle
AW\big(2,\,2,\,0,\,0,\,1,\,0,\,\,-4u(u+d+1)\,v,\,-(u+1)(u+d)-2v^2-\frac{d^2}2\big).
\end{eqnarray}
\item For the Krawtchouk case:
\begin{equation}
AW\big(2,\,0,\,0,\,1,\,1,\,2v-1,\,0,\,0\big).
\end{equation}
\item For the Bannai-Ito case, if $d$ is even:
\begin{equation}
AW\big(-2,\,0,\,0,\,1,\,1,\,4uu^*\!-2(d\!+\!1)\,v,\,2uv-(d\!+\!1)\,u^*,\,2u^*v-(d\!+\!1)\,u\big).
\end{equation}
\item For the Bannai-Ito case, if $d$ is odd:
\begin{eqnarray}
&&\hspace{-24pt}\textstyle AW\big(-2,\,0,\,0,\,1,\,1,\,
-2u^2-2u^*{}^2+2v^2+\frac{(d+1)^2}2, \nonumber \\
&&\textstyle
-u^2+u^*{}^2-v^2+\frac{(d+1)^2}4,\,u^2\!-u^*{}^2\!-v^2\!+\frac{(d+1)^2}4\big).
\end{eqnarray}
\end{itemize}
\end{lemma}
\proof Direct computations with formulas (\ref{betap1})--(\ref{asipq2}).\qed

\begin{remark} \label{orphan} \rm 
In the characteristic 2, the normalization of part 1 in Lemma
\ref{normrules} is not available, hence our results are incomplete if char
$\fK=2$. In particular, we miss the orphan case (with char $\fK=2$ and $d=3$)
completely. A general parameter array of the orphan type  
can be normalized as follows:
\begin{eqnarray*}
\left(\theta_0,\theta_1,\theta_2,\theta_3\right)\equal\left(0,\,s+1,\,1,\,s\right),\\
\left(\theta^*_0,\theta^*_1,\theta^*_2,\theta^*_3\right)\equal\left(0,\,s^*\!+1,\,1,\,s^*\right),\\
\left(\varphi_1,\varphi_2,\varphi_3\right)\equal\left(r,\,1,\,r+s+s^*\right),\\
\left(\phi_1,\phi_2,\phi_3\right)\equal\left(r+s+s\,s^*,\,1,\,r+s^*\!+s\,s^*\right).
\end{eqnarray*}
Here adjusted $\theta_0=0$, $\theta^*_0=0$, $h=1$, $h^*=1$ in \cite[Example
35.14]{terwgen}. The Askey-Wilson relations are
\begin{equation} \label{orphanaw}
AW(0,\,1,\,1,\,s^2+s,\,s^*{}^2+s^*,\,s\,s^*,\,r\,s,\,r\,s^*).
\end{equation}
The relations can be renormalized to $\rho=0$, $\rho^*=0$ by affine translations (in 4
ways, generally). The normalized coefficients $\eta,\eta^*,\omega$ in
(\ref{orphanaw}) are dependent on two free parameters, so there is a
relation between them. Here is the relation, in the form invariant under
affine rescaling:
\begin{eqnarray} \label{orphanwc}
\left(\omega^2\!-\varrho\varrho^*\right)^2=\omega\,(\gamma\omega-\gamma^*\varrho)
(\gamma^*\omega-\gamma\varrho^*).
\end{eqnarray}
\end{remark}

\section{Classification of AW relations}
\label{classifics}

Askey-Wilson relations can be consistently classified by families of
orthogonal polynomials in the same way as Leonard pairs. The classification
is presented in the first two columns of Table \ref{elltab}.
\begin{table}
\begin{center} \begin{tabular}{|c|c|c|}
\hline \centering Askey-Wilson type & Askey-Wilson coefficients & with Leonard pairs \\ \hline %
$q$-Racah & $\beta\neq\pm2$, $\underline{\gamma=\gamma^*\!=0}$,
 $\widehat{\varrho}\,\widehat{\varrho}^*\!\neq 0$ & --- \\
$q$-Hahn & $\beta\neq\pm2$, $\underline{\gamma=\gamma^*\!=0}$,
 $\widehat{\varrho}=0$, $\widehat{\varrho}^*\widehat{\eta}\neq 0$ & --- \\
Dual $q$-Hahn & $\beta\neq\pm2$, $\underline{\gamma=\gamma^*\!=0}$,
 $\widehat{\varrho}^*\!=0$, $\widehat{\varrho}\,\widehat{\eta}^*\!\neq 0$ & --- \\
 $q$-Krawtchouk & $\beta\neq\pm2$, $\underline{\gamma=\gamma^*\!=0}$,
 $\widehat{\varrho}=\widehat{\eta}=0$ & $\widehat{\varrho}^*\widehat{\eta}^*\neq 0$\\
Dual $q$-Krawtchouk & $\beta\neq\pm2$, $\underline{\gamma=\gamma^*\!=0}$,
 $\widehat{\varrho}^*\!=\widehat{\eta}^*\!=0$ & $\widehat{\varrho}\;\widehat{\eta}\neq 0$ \\
$\begin{array}{l} \mbox{Quantum/affine}\vspace{-5pt}\\
 \mbox{$q$-Krawtchouk}\end{array}$ & $\beta\neq\pm2$, $\underline{\gamma=\gamma^*\!=0}$,
 $\widehat{\varrho}=\widehat{\varrho}^*\!=0$ & $\widehat{\eta}\,\widehat{\eta}^*\!\neq 0$ \\
Racah & $\beta=2$, $\gamma\,\gamma^*\neq 0$, $\underline{\varrho=\varrho^*=0}$ & --- \\
Hahn & $\beta=2$, $\gamma=0$, $\gamma^*\neq 0$,
 $\underline{\varrho^*=0}$ & $\varrho\neq 0$, $\underline{\omega=0}$\\
Dual Hahn & $\beta=2$, $\gamma^*\!=0$, $\gamma\!\neq 0$,
 $\underline{\varrho=0}$ & $\varrho^*\!\neq 0$, $\underline{\omega=0}$\\
Krawtchouk & $\beta=2$, $\gamma=\gamma^*\!=0$ & $\varrho\varrho^*\!\neq 0$,
$\underline{\eta=\eta^*\!=0}$ \\
Bannai-Ito & $\beta=-2$, $\underline{\gamma=\gamma^*=0}$ &
 $\widehat{\varrho}\,\widehat{\varrho}^*\neq 0$ \vspace{1pt}\\
\hline
\end{tabular} \end{center}
\caption{Classification of Askey-Wilson relations} \label{elltab}
\end{table}
In each line, the underlined equalities can be achieved by using affine
translations if the preceding conditions are satisfied. If $\beta\neq 2$,
by $\widehat{\varrho}$, $\widehat{\varrho}^*$, $\widehat{\omega}$,
$\widehat{\eta}$, $\widehat{\eta}^*$ we denote the Askey-Wilson coefficients
in a normalization specified by part 1 of Lemma \ref{normrules}.

The first part of the following theorem establishes the consistency of
Askey-Wilson types for Leonard pairs and for Askey-Wilson relations.
\begin{theorem} \label{classth}
Assume that $d\ge 3$. Let $(A,A^*)$ denote a Leonard pair on
$V$, and let $AW$ denote the Askey-Wilson relations satisfied by $(A,A^*)$.
\begin{enumerate}
\item The Askey-Wilson relations $AW$ have the same Askey-Wilson type as
the Leonard pair $(A,A^*)$.%
\item If there is other Leonard pair on $V$ that satisfies AW, it has the
same Askey-Wilson type as $(A,A^*)$.%
\item There exist unique affine translation which normalizes $AW$ according
to the specifications of Lemma $\ref{normrules}$. 
\item The Askey-Wilson relations AW satisfy all inequalities in the last two
columns of Table $\ref{elltab}$ on the corresponding line. All underlined
equalities can be achieved after an affine translation, and such an affine
translation is unique. The indicated nonzero coefficients can be normalized
to any chosen values by an affine scaling.
\end{enumerate}
\end{theorem}
\proof For the first statement, check the results in Section
\ref{qpararrays1} (or Section \ref{qpararrays2}) and Section
\ref{otherpararrays}, and observe that the Askey-Wilson relations associated
to any parameter array have the same Askey-Wilson type as the parameter
array, with the exception of the ambiguity between the quantum
$q$-Krawtchouk and affine $q$-Krawtchouk types.

The second statement is an immediate consequence.

For the third statement, we have to prove that cases $6$ and $7$ of Lemma
$\ref{normrules}$ do not apply to $AW$. Assuming the contrary, $AW$ would
have the Krawtchouk type. In the corresponding normalized form of Lemma
\ref{awnormals2} we would have $v\in\{0,1\}$. But then the Krawtchouk
parameter array of Lemma \ref{normlps2} degenerates, since $\phi_i=0$ or
$\psi_i=0$ for all $i=1,2,\ldots,d$. The third statement follows.

The inequalities of the last column of Table \ref{elltab} can be checked by
inspecting all Askey-Wilson relations in Lemmas \ref{awnormals},
\ref{awnormals1} and \ref{awnormals2}. Normalization by affine translations
follows from the Lemma \ref{normrules} and the previous part here.
Normalization by affine scalings is clear.\qed\\

Note that the normalization specified by Lemma \ref{normrules} follows
implications of part 4 of Theorem \ref{classth}. By part 3 of Theorem
\ref{classth}, there is a unique affine translation to set two specified
Askey-Wilson coefficients to zero. For each type of Leonard pairs, we get
two Askey-Wilson coefficients which are certainly nonzero after the
normalizing affine translation. These coefficients can be characterized as
follows: they are the first nonzero (after the normalizing translation)
coefficients in the two sequences
\begin{equation} \label{scalingpars}
(\gamma,\,\varrho,\,\eta,\,\eta^*)\quad\mbox{and}\quad
(\gamma^*,\varrho^*,\eta^*,\eta).
\end{equation}
By affine scalings, the two coefficients can be normalized to any convenient
values. In the Askey-Wilson relations of Lemmas \ref{awnormals1} and
\ref{awnormals2}, the normalized values depend only on $\beta$:
\begin{eqnarray} \label{scnormal}
\mbox{ }\gamma,\gamma^*: && \mbox{$2$ (if $\beta=2$)};\nonumber\\
\mbox{ }\varrho,\varrho^*: && \left\{\begin{array}{cl}
4\!-\!\beta^2, & \mbox{if } \beta\neq\pm 2,\\
1, & \mbox{if } \beta=\pm 2;\end{array}\right.\\
\mbox{ }\eta,\eta^*: &&
\left\{\begin{array}{cl}\sqrt{\beta\!+\!2}\,(\beta\!-\!2), & \mbox{if }
\eta\eta^*\neq 0 \mbox{ or } \omega=0,\\
\sqrt{\beta\!+\!2}\,(\beta\!-\!2)\,Q_{d+1}, & \mbox{if } \eta\eta^*=0 \mbox{
and } \omega\neq 0.\end{array}\right.\nonumber
\end{eqnarray}
$Q_{d+1}$ can be independently defined by the linear recurrence
$Q_{n+2}=\beta\,Q_n-Q_{n-2}$ with the initial values
$Q_{-1}=Q_1=\sqrt{\beta\!+\!2}$, $Q_0=2$, $Q_2=\beta$. One can take for
$\sqrt{\beta\!+\!2}$ any of the two values of the square root. In the context of
Lemma \ref{awnormals1}, we should identify $\sqrt{\beta+2}$ with $q+q^{-1}$.
The effect of changing the sign of $\sqrt{\beta\!+\!2}$ is
multiplication of $A$ and/or $A^*$ by $-1$. 

\section{Uniqueness of normalizations}
\label{conclusions}

The results in Sections \ref{qpararrays1} through \ref{otherpararrays} can be
used to compute the Askey-Wilson relations for any Leonard pair. 
To do this, one may take a parameter array corresponding to a given
Leonard pair; then find an affine transformation (\ref{afftr}) which normalizes
the parameter array by (\ref{scalingtr}) to one of the forms of Lemmas
\ref{normlps}, \ref{normlps1} or \ref{normlps2}; then pick up the corresponding normalized
Askey-Wilson relations in Lemmas \ref{awnormals}, \ref{awnormals1} or \ref{awnormals2};
and then apply the inverse affine transformation to the normalized relations using formula
(\ref{newaw}). This procedure can be applied for any $d$, although for $d<3$
the type of a representing parameter array is ambiguous and the Askey-Wilson
relations are not unique.

For the rest of this section, we refer to the results of Sections
\ref{qpararrays2} and \ref{otherpararrays}. We assume $d\le 3$ and adopt the
following terminology. A pair of Askey-Wilson relations is called {\em
normalized} if it satisfies the specifications of Lemma \ref{normrules} and
the description in the previous section; see (\ref{scalingpars}) and
(\ref{scnormal}). A Leonard pair is {\em normalized} if it satisfies
normalized Askey-Wilson relations. A parameter array is {\em normalized} if
it can be expressed in one of the forms of Lemma \ref{normlps1} or Lemma
\ref{normlps2}.

We consider the following questions:
\begin{question}\rm \label{question0}
How unique is normalization of Askey-Wilson relations?
\end{question}
\begin{question}\rm \label{question1}
Given a Leonard pair, how unique is its normalization?
\end{question}
\begin{question} \rm \label{question4}
Is every normalized Leonard pair representable by a normalized parameter
array?
\end{question}
\begin{question}\rm \label{question2}
Are normalized parameter arrays represented uniquely by the
forms in Lemmas \ref{awnormals1} and \ref{awnormals2}?%
\end{question}
\begin{question} \rm \label{question3}
Do the relation operators $\downarrow$, $\Downarrow$,
$\downarrow\Downarrow$ preserve the set of normalized parameter arrays? %
\end{question}
\begin{table}
\begin{center} \begin{tabular}{|c|c|c|}
\hline Askey-Wilson & Affine scaling & Conversion of normalized \vspace{-3pt}\\
type & $(t,t^*)$ &  parameter arrays \\ \hline
$q$-Racah 
 & $(-1,1)$  & $s\mapsto -s$, $r\mapsto -r$\vspace{-2pt}\\
 & $(1,-1)$  & $s^*\mapsto -s^*$, $r\mapsto -r$\\
$q$-Hahn & $(\sqrt{-1},-1)$  & $s^*\mapsto-s^*$, $r\mapsto\sqrt{-1}\,r$\\
Dual $q$-Hahn & $(-1,\sqrt{-1})$  & $s\mapsto-s$, $r\mapsto\sqrt{-1}\,r$\\
$q$-Krawtchouk & $(1,-1)$ & $s^*\mapsto -s^*$\\
Dual $q$-Krawtchouk & $(-1,1)$ & $s\mapsto -s$\\
$\begin{array}{l} \mbox{Quantum and affine}\vspace{-5pt}\\
 \mbox{$q$-Krawtchouk}\end{array}$ & $(\zeta_3,\zeta_3)$ &  $r\mapsto\zeta_3\,r$\\
Racah & --- & --- \\
Hahn & $(-1,1)$ & $\Downarrow$ and $v\mapsto-v$\\
Dual Hahn & $(1,-1)$ & $\downarrow$ and $v\mapsto-v$\\
Krawtchouk & $(-1,1)$ & $\Downarrow$ \vspace{-2pt}\\
 & $(1,-1)$ & $\downarrow$\\
Bannai-Ito
 & $(-1,1)$ & If $d$ even: $\Downarrow$ and $u\mapsto-u$, $v\mapsto-v\;\;$\vspace{-2pt}\\
 & $(1,-1)$ & If $d$ even: $\downarrow$ and $u^*\!\mapsto-u^*$, $v\mapsto-v$\\
\hline
\end{tabular} \end{center}
\caption{Reparametrization of different normalizations} \label{normtr}
\end{table}
\begin{table}
\begin{center} \begin{tabular}{|c|c|c|}
\hline Askey-Wilson & Change of sign  &  Parameter array \vspace{-3pt} \\
 type & of $\sqrt{\beta\!+\!2}$ & stays invariant  \\ \hline %
$q$-Racah & --- & $r\mapsto 1/r$; also (\ref{qracahinv}), (\ref{qracahin2}) \\
$q$-Hahn & $q\mapsto -q$, $s^*\mapsto (-1)^{d+1}s^*$ & ---  \\
Dual $q$-Hahn & $q\mapsto-q$, $s\mapsto (-1)^{d+1}s$ & --- \\
$q$-Krawtchouk & If $d$ odd: $q\mapsto-q$, $s^*\mapsto-s^*\!$  & If $d$ even: $q\mapsto -q$\\
Dual $q$-Krawtchouk & If $d$ odd: $q\mapsto -q$, $s\mapsto-s\;\;$ & If $d$ even: $q\mapsto-q$ \\
$\begin{array}{l} \mbox{Quantum and affine}\vspace{-5pt}\\
 \mbox{$q$-Krawtchouk}\end{array}$ & $q\mapsto-q$, $r\mapsto (-1)^{d+1}r$ & --- \\
Racah & --- & $v\mapsto -v-d-1$ \\
Bannai-Ito & --- & If $d$ odd: $v\mapsto -v$ \\
\hline
\end{tabular} \end{center}
\caption{Alternative normalization and invariant reparametrizations}
\label{invartr}
\end{table}
Regarding the first question, non-uniqueness occurs for two reasons:
\begin{itemize}
\item There exist affine scalings by small roots of unity that leave the
first nonzero coefficients in both sequences of (\ref{scalingpars}) invariant.
The list of these affine scalings is given by the first two columns of Table \ref{normtr}. 
By $\zeta_3$ we denote a primitive cubic root of unity. In the
$q$-Racah, Krawtchouk and Bannai-Ito cases, two given scalings can be
composed. In the $q$-Hahn, dual $q$-Hahn and quantum/affine $q$-Krawtchouk
cases, there are non-trivial iterations of the given scalings. 
The third column of Table \ref{normtr} gives corresponding conversions
of parameter arrays.
\item In  all $q$-Hahn and $q$-Krawtchouk cases, there exists an alternative normalization
of the two nonzero Askey-Wilson coefficients from (\ref{scalingpars}), 
with the other sign of $\sqrt{\beta\!+\!2}$. This effectively multiplies $A$ or $A^*$ (or both) by $-1$.
The corresponding action on parameter arrays is given by the second  column of Table \ref{invartr}.
\end{itemize}
Normalization of Askey-Wilson relations is unique in the Racah case. Otherwise, the normalization is
unique in the $q$-Racah, Hahn, dual Hahn, Krawtchouk and Bannai-Ito cases
if (and only if) the Askey-Wilson relations remain invariant under the respective affine scalings.
This means that the non-scaled coefficients (such as $\widehat{\omega}$, $\widehat{\eta}$,
$\widehat{\eta}^*$ in the $q$-Racah and Bannai-Ito cases) in the normalized relations are equal to zero.

Question \ref{question1} is equivalent to Question \ref{question0}. However, existence and
uniqueness of representation of a normalized Leonard pair by a normalized parameter array is determined by Questions \ref{question4} through \ref{question3}. 

An important discrepancy between
normalization of Askey-Wilson relations and normalization of parameter arrays occurs in
the Bannai-Ito case if $d$ is odd. As Table \ref{normtr} implies, affine scalings that preserve
normalization of Askey-Wilson relations cannot be realized by transformations of normalized
parameter arrays then. This has implications for Question \ref{question4}.
\begin{lemma} \label{normscale}
\begin{enumerate}
\item Any normalized Leonard pair can be represented by a normalized
parameter array, 
except when the Askey-Wilson type is Bannai-Ito, and $d$ is odd.%
\item Suppose that $d$ is odd. Let $(B,B^*)$ denote the Leonard pair
represented by the parameter array of the Bannai-Ito type in Lemma $\ref{normlps2}$.
Then the following four Leonard pairs satisfy 
normalized Askey-Wilson relations of the Bannai-Ito type:
\begin{equation} \label{bito4}
(B,B^*),\quad (-B,B^*),\quad (B,-B^*),\quad (-B,-B^*).
\end{equation}
Of these Leonard pairs, only $(B,B^*)$ can be represented by a normalized
parameter array.
\end{enumerate}
\end{lemma}
\proof Let $(A,A^*)$ denote a normalized Leonard pair on $V$. Let $\Phi$
denote a parameter array for $(A,A^*)$. Let $\Phi^\#$ denote a normalization
of $\Phi$ by (\ref{scalingtr}); it can be expressed in one of the forms of
Lemmas \ref{awnormals1} and \ref{awnormals2}. The parameter arrays $\Phi$
and $\Phi^\#$ differ by an affine scaling from Table \ref{normtr}, plus (in
some $q$-cases) possibly the change of the sign of $\sqrt{\beta\!+\!2}$ in
the Askey-Wilson relations. If the sign of $\sqrt{\beta\!+\!2}$ is changed,
one can apply a corresponding reparametrization in the second column of Table \ref{invartr}. 
Reparametrizations for relevant affine scalings are indicated in Table \ref{normtr},
except for the Bannai-Ito case with odd $d$. Hence $\Phi$ is normalized as
well, except perhaps when it has the Bannai-Ito type and $d$ is odd.

Now we prove the second part. The four Leonard pairs in (\ref{bito4}) are normalized
according to our discussion of Questions \ref{question0} and \ref{question1}.
The Bannai-Ito parameter array of Lemma $\ref{normlps2}$
has the following property: the even-indexed $\theta_i$'s and the even indexed
$\theta^*_i$'s form increasing sequences. Since $d$ is assumed odd, 
the relation operators $\downarrow$, $\Downarrow$,
$\downarrow\Downarrow$ preserve this property. But affine scalings
(\ref{scalingtr}) with $t=-1$ or $t^*=-1$ reverse this property
for $\theta_i$'s or $\theta^*_i$'s, respectively. Hence, in all parameter arrays representing
$(-B,B^*)$, $(B,-B^*)$ or $(-B,-B^*)$ the even-indexed $\theta_i$'s and/or
the even indexed $\theta_i^*$'s are in the decreasing order. The conclusion is that
the four Leonard pairs $(\pm B,\pm B^*)$ cannot be transformed
to each other by change of the parameters $u,u^*,v$ or the relation operators.  
Hence only $(B,B^*)$ can be represented as a specialization of the
Bannai-Ito parameter array  of Lemma $\ref{normlps2}$.\qed

Questions \ref{question2} and \ref{question3} determine how unique are
representations of normalized Leonard pairs by normalized parameter arrays.
Invariant reparametrization of parameter arrays do occur. They are given in
the third column of Table \ref{invartr}. In the $q$-Racah case, we additionally have
the following invariant transformations:
\begin{eqnarray} \label{qracahinv}
&q\mapsto 1/q, \qquad s\mapsto 1/s, \qquad s^*\mapsto1/s^*;\\
\label{qracahin2} &q\mapsto -q, \quad s\mapsto (-1)^ds, \quad
s^*\mapsto(-1)^ds^*,\quad r\mapsto (-1)^{d+1}r.
\end{eqnarray}
Question \ref{question3} is thoroughly answered in Table \ref{invartr2}.
There, ``Switch" means interchanging the quantum $q$-Krawtchouk and affine
$q$-Krawtchouk types of parameter arrays. Of course, $\downarrow\Downarrow$
is the composition of $\downarrow$ and $\Downarrow$. As wee see, the
relation operators preserve normalization of parameter arrays in all
$q$-cases, in the Racah case, and in the Bannai-Ito case with odd $d$.
\begin{table}
\begin{center} \begin{tabular}{|c|c|c|}
\hline Askey-Wilson type & Conversion to $\Downarrow$ & Conversion to $\downarrow$ \\ \hline %
$q$-Racah &  $s\mapsto 1/s$ & $s^*\mapsto 1/s^*$ \\
$q$-Hahn & $q\mapsto 1/q$, $s^*\!\mapsto 1/s^*$ & $s^*\mapsto 1/s^*$ \\
Dual $q$-Hahn & $s\mapsto 1/s$ & $q\mapsto 1/q$, $s\mapsto 1/s$ \\
$q$-Krawtchouk & $q\mapsto 1/q$, $s^*\!\mapsto 1/s^*$ & $s^*\mapsto 1/s^*$ \\
Dual $q$-Krawtchouk & $s\mapsto 1/s$ & $q\mapsto 1/q$, $s\mapsto 1/s$  \\
Quantum/affine $q$-Krawtchouk & Switch & Switch and $q\mapsto1/q$\\
Racah & $u\mapsto -u-d-1$ & $u^*\mapsto -u^*\!-d-1$ \\
Hahn & --- & $u^*\mapsto -u^*\!-d-1$ \\
Dual Hahn &  $u\mapsto -u-d-1$ & --- \\
Bannai-Ito, $d$ odd & $u\mapsto-u$ & $u^*\mapsto -u^*$ \\
\hline
\end{tabular} \end{center}
\caption{Relative parameter arrays} \label{invartr2}
\end{table}

\bibliographystyle{alpha}
\bibliography{../terwilliger,../../hypergeometric}

\end{document}